\documentclass[12pt]{article}
\usepackage[affil-it]{authblk}
\usepackage{amsmath}
\usepackage{chngcntr}
\counterwithin{equation}{section}
\usepackage{times}
\usepackage{mathptmx}
\title{SOLUTION PROPERTIES  FOR PERTURBED LINEAR AND NONLINEAR INTEGRAL EQUATIONS}
\author{M.Z. TSOUKALAS$^1$, P.G. ASTERIS$^{1,*}$\\(martsoukalas@ta.aspete.gr,panagiotisasteris@gmail.com)}
\affil{Computational Mechanics Laboratory, School of Pedagogical and Technological Education, Athens, Greece      \\
$^*$ Corresponding author         }

\date{\today}
\begin{document}
\newtheorem{theo}{Theorem}[section]
\newtheorem{zlemma}{Lemma}[section]
\newtheorem{zlem}{Lemma}[section]
\newtheorem{zdef}{Definition}[section]
\newtheorem{zcorr}{Corrolary}[section]
\newtheorem{prop}{Proposition}[section]
\newtheorem{zexam}{Example}[section]
\maketitle
\begin{abstract}
In this study we consider  perturbative series solution with respect to a parameter $\epsilon>0$. In this methodology the solution is considered as an infinite sum of a series of functional terms which usually converges fast to the exact desired solution. Then we investigate perturbative solutions for kernel perturbed integral equations and prove the convergence in an appropriate ranges of the perturbation series. Next  we investigate perturbation series solutions for nonlinear perturbations of integral equations of Hammerstein type and formulate conditions for their convergence. Finally we prove the existence of a maximal perturbation range for  non linear integral equations. \end{abstract}
\pagebreak
\section{\normalsize \bf Introduction}\label{sec0}
Integral equations arise in many fields in mathematical physics, biology, chemical kinetics, mechanics, etc.. In recent years there is a literature dealing with homotopy perturbation methods in integral equations (\cite{agar,ahakhi,ahsa,alma,alma2,chamm2,esku, einteg2,gsinteg1,ihamm1,nawa,yinteg3}). However despite its fast numerical convergence, the convergence of the  OHAM (optimal homotopy asymptotic method) is not proved formally in (\cite{agar,ahsa,alma,alma2,ahakhi,einteg2,esku,nawa,yinteg3}). In this work we deal mainly with existence theorems for nonlinear integral equations for solutions in $C[0,1]$ or $L^2[0,1]$ spaces. Existence theorems are given in (\cite{hoch,precu,regan}) but here we give versions of these theorems for perturbed integral equations by formulating and proving basic existence theorems for non linear integral equations considering perturbative series solution with respect to a parameter $\epsilon>0$. In this methodology the solution is considered as an infinite sum of a series of functional terms which usually converges fast to the exact desired solution. Particularly we investigate perturbative solutions for kernel perturbed integral equations and prove the convergence in  appropriate ranges of the perturbation series. In the next section \ref{sec3} we investigate perturbation series solutions for nonlinear perturbations of non linear integral equations of Hammerstein type and formulate conditions for their convergence. Finally we prove the existence of a maximal perturbation range for  non linear integral equations. The methods used here can be exploited to implement numerical procedures for perturbation series solutions.\\\\
\section{Series Solutions for perturbed kernels}\label{sec2}
Here we consider 'small' in a sense as we shall see perturbations in a non linear sense of the kernel $K$ of the integral equation and we prove that the perturbative series solution converges absolutely in an appropriate range for $\epsilon$, the perturbation parameter.\\
In theorem (\ref{pertheo1}) we deal with perturbations of the kernel in the linear Fredholm equation of the second kind and prove that if the unperturbed original equation has a $C[0,1]$ or $L^2[0,1]$ solution then the perturbed equation has a corresponding $C[0,1]$ or $L^2[0,1]$ solution. Next in theorem (\ref{pertheo2a}) we prove the corresponding result dropping the boundedness assumption for the perturbed kernel. Continuing we prove in theorem (\ref{pertheo2}) that the corresponding result holds under suitable assumptions for a non linear  Hammerstein integral equation.  Following we prove the corresponding result in theorem (\ref{pertheo1dd}) for a linear Fredholm equation of the second kind. with boundedly differentiable kernel. Next we prove the theorem (\ref{pertheo2d2}) stating the result for $L^2$ integrable perturbed kernel in a Hammerstein equation. We prove the same kind of result for equations with boundedly differentiable kernel in theorem (\ref{pertheo2d}). We proceed with theorem (\ref{pertheo3d}) stating the corresponding result for $L^1$ integrable kernel. We denote by $\phi(\epsilon,x)=\phi_{(0,0)}(x)$ the unknown sought function and by $\phi_{(0,\nu)}(x)=\frac{\partial^{\nu}\phi(\epsilon,x)}{\partial \epsilon^{\nu}}|_{\epsilon=0},\nu=0,1,2,...$ its various order derivatives at $\epsilon=0$.\\\\
{\theo{Let's consider the perturbed integral equation
\begin{eqnarray}
&&  \phi(\epsilon,x)-\omega\int_0^1\Gamma (\epsilon,x,y)\phi(\epsilon,y)dy=f(x),\;\;\;\;x\in[0,1],\epsilon\geq 0\label{eqinte1nla}
                                      \end{eqnarray}
                                      where $\Gamma (\epsilon,x,y)=\Gamma_0(x,y)+\epsilon \Gamma_1(x,y)$ is continuous and $L^2$ integrable respectively functions of their variables and $\int_0^1\int_0^1|\Gamma (\epsilon,x,y|dxdy<\infty$.
                                       Suppose that
\begin{eqnarray}
 && |\Gamma_0(\epsilon,x,y)|\leq C,\;\;\;\;\;\;0\leq x,y\leq 1,\epsilon\geq 0.
\end{eqnarray}
Assuming that the unperturbed equation for $\epsilon=0$ has a solution in $C([0,1])$ or $L^2([0,1])$ respectively
then the perturbation series
\begin{eqnarray}
&&\phi_{(0,0)}(x)=\omega \int_0^1\Gamma_0,x,y)\phi(0,y)dy+f(x)\label{per0}\\
&&\phi_{(0,1)}(x)=\omega \int_0^1\Gamma_1(x,y)\phi(0,y)dy+\omega \int_0^1\Gamma_0(x,y)\phi_{(0,1)}(y)dy\label{per1}\\
&&\cdot\cdot\cdot\nonumber\\
&&\phi_{(0,\nu)}(x)=\frac{\partial^{\nu}\phi(\epsilon,x)}{\partial \epsilon^{\nu}}|_{\epsilon=0}=\omega \int_0^1\Gamma_1((x,y)\phi_{(0,\nu-1)}(y)dy\nonumber\\
&&+\omega \int_0^1\Gamma (0,x,y)\phi_{(0,\nu)}(y)dy,\nu=1,2,....\label{pern}
\end{eqnarray}
converges absolutely  in  $||\cdot(x)||_0=\sup_{x\in[0,1]}|\cdot(x)|$ norm and $||\cdot||_2$ norm to a solution
\begin{eqnarray}
&&\phi(\epsilon,x)=\sum_{j=0}^{\infty}\phi_{(0,j)}(x)\frac{\epsilon^j}{j!}\label{serpera111}
\end{eqnarray}
which is continuous with respect to $(\omega,\epsilon,x)$.}\label{pertheo1}\\\\}
{\bf Proof:} As $\Gamma_0(x,y)\leq C$ and $\Gamma_1(x,y)\leq C$ we have that $||\Gamma_0||_1\leq C, ||\Gamma_1||_1\leq C$, which implies
\begin{eqnarray}
&&||\phi_0^{(\nu)}||_0=||\phi(0,x)||_0\leq |\omega| ||\Gamma_1||_1\cdot ||\phi_0^{(\nu-1)}||_0+|\omega| ||\Gamma_0||_1\cdot ||\phi_0^{(\nu)}||_0\label{ineqper01}
\end{eqnarray}
leading to
\begin{eqnarray}
&&||\phi_0^{(\nu)}||_0\leq\frac{|\omega|||\Gamma_1||_1}{1-|\omega|||\Gamma_0||_1}||\phi^{(\nu-1)}||_2\label{ineqper02}
\end{eqnarray}
which implies that for $$\frac{|\omega|||\Gamma_1||_1}{1-|\omega|||\Gamma_0||_1}=\rho_0,$$
we have that
\begin{eqnarray}
&&||\phi(\epsilon,x)||_0=\sum_{j=0}^{\infty}||\phi_{(0,j)}(x)||_0\frac{\epsilon^j}{j!}\\
&&\leq \sum_{j=0}^{\infty}\rho_0^{j}\frac{\epsilon^j}{j!}||\phi_0||_0\label{serper02}
\end{eqnarray}
the solution converges in $||\cdot||_0$ for any value of $\omega,\rho_0,\epsilon\geq 0$.\\\\
Also we observe that as $||\Gamma_j||_2<\infty$ for $j=0,1$,
\begin{eqnarray}
||\phi_0^{(\nu)}||_2&=&||\phi_{(0,0)}(x)||_2\leq |\omega| ||\Gamma_1||_2\cdot ||\phi_{(0,\nu-1)}||_2\nonumber\\
&+&|\omega| ||\Gamma_0||_2\cdot ||\phi_{(0,\nu)}||_2\label{ineqper1}
\end{eqnarray}
leading to
\begin{eqnarray}
&&||\phi_{(0,\nu)}||_2\leq\frac{|\omega|||\Gamma_1||_2}{1-|\omega|||\Gamma_0||_2}||\phi_{(0,\nu-1)}||_2\label{ineqper2}
\end{eqnarray}
which implies that for $$\frac{|\omega|||\Gamma_1||_2}{1-|\omega|||\Gamma_0||_2}=\rho,$$
we have that
\begin{eqnarray}
&&||\phi(\epsilon,x)||_2\leq\sum_{j=0}^{\infty}||\phi_{(0^,j)}(x)||_2\left[\frac{\epsilon^j}{j!}\right]\\
&&\leq \sum_{j=0}^{\infty}|\rho^{j}\left[\frac{\epsilon^j}{j!}\right]||\phi_0||_2\leq ||\phi_0||_2e^{|\rho \epsilon|}\label{serper2}
\end{eqnarray}
the solution converges for any value of $\omega,\rho,\epsilon\geq 0$.\\\\
As a convergent power series everywhere it is continuous with respect to $(x,\epsilon)$ as one can see immediately from uniform convergence.
Also we have that is continuous with respect to $\omega$.\\\\\\
{\theo{Let's consider the perturbed integral equation
\begin{eqnarray}
&&  \phi_{(\epsilon)}(x)-\omega\int_0^1\Gamma (\epsilon,x,y)\phi_{(\epsilon)}(y)dy=f(x),\;\;\;\;x\in[0,1],\epsilon\geq 0\label{eqinte1nl}
                                      \end{eqnarray}
                                      where $\Gamma (\epsilon,x,y)=\Gamma_0(x,y)+\epsilon \Gamma_1(x,y)$ is  $L^2$ integrable respectively functions of their variables and $\int_0^1\int_0^1|\Gamma (\epsilon,x,y)|dx dy<\infty$.
%
Assuming that the unperturbed equation for $\epsilon=0$ has a solution in $C([0,1])$
then the perturbation series
\begin{eqnarray}
&&\phi_{(0,0)}(x)=\omega \int_0^1\Gamma_0,x,y)\phi_{(0,0)}(y)dy+f(x)\label{per1a}\\
&&\phi_{(0,1)}(x)=\omega \int_0^1\Gamma_1(x,y)\phi_{(0,0)}(y)dy+\omega \int_0^1\Gamma_0(x,y)\phi_{(0,1)}(y)dy\label{per1a}\\
&&\cdot\cdot\cdot\nonumber\\
&&\phi_{(0,\nu)}(x)=\frac{\partial^{\nu}\phi(x)}{\partial x^{\nu}}|_{x=0}=\omega \int_0^1\Gamma_1((x,y)\phi_{(0,\nu-1)}(y)dy\nonumber\\
&&+\omega \int_0^1\Gamma (0,x,y)\phi_{(0,\nu)}(y)dy,\nu=1,2,....\label{perna}
\end{eqnarray}
%
converges absolutely  in      $||\cdot||_0=\sup_{x\in[0,1]}|\cdot(x)|$ norm and $||\cdot||_2$ norm to a solution
\begin{eqnarray}
&&\phi(\epsilon,x)=\sum_{j=0}^{\infty}\phi_{(0,j)}(x)\frac{\epsilon^j}{j!}\label{serpera}
\end{eqnarray}
which is continuous with respect to $(\omega,\epsilon,x)$.}\label{pertheo2a}\\\\}
{\bf Proof:}  As a result of the $L^2$ integrability we have that we have converging in $L^2$ sequences of kernels to $\Gamma_0,\Gamma_1$ respectively such that  $\Gamma_{0,n}(x,y)\leq C_n\rightarrow C_{|infty}\leq\infty$ and $\Gamma_{1,n}(x,y)\leq C_n\rightarrow C_{\infty,1}\leq\infty$ and we have that $||\Gamma_{0,n}||_1\leq C_n, ||\Gamma_{1,n}||_1\leq C_n$, which implies corresponding inequalities to (\ref{ineqper01},\ref{ineqper02}), and therefore there exists a solution sequence
\begin{eqnarray}
&&||\phi_{(n,\epsilon)}(x)||_0=\sum_{j=0}^{\infty}||\phi_{(n,0,j)}(x)||_0\frac{\epsilon^j}{j!}\\
&&\leq \sum_{j=0}^{\infty}\rho_{0,n}^{j}\frac{\epsilon^j}{j!}||\phi_{0,n}||_0,n=0,1,2,...\label{serper0n2aa}
\end{eqnarray}
the solution converges in $||\cdot||_0$ for any value of $\omega,\rho_0,\epsilon\geq 0$.\\\\
As $\rho_{0,n}\rightarrow \rho_0<\infty$ we have that the limit solution $||\phi(\epsilon,x)||_0<\infty$ is bounded  everywhere for $(\omega,x,\epsilon)$.
As the convergence is uniform we have that $\phi(\epsilon,x)$ is continuous.\\\\
In the same manner we have as a result that $\phi(\epsilon,x)$ belongs to $L_2[0,1]$.\\\\
{\theo{Let's consider the perturbed integral equation
\begin{eqnarray}
&&  \phi(\epsilon,x)-\omega\int_0^1\Gamma (\epsilon,x,y)\psi(y,\phi(\epsilon,y))dy=f(x),\;\;\;\;x\in[0,1],\epsilon\geq 0\label{eqinte1nol}
                                      \end{eqnarray}
                                      where $\Gamma (\epsilon,x,y)=\Gamma_0(x,y)+\epsilon \Gamma_1(x,y)$ is  $L^2$ integrable respectively function of its variables and $\int_0^1\int_0^1|\Gamma (\epsilon,x,y)|^2dx dy<\infty$.
                                       Suppose that
\begin{eqnarray}
 && |\Gamma_j(\epsilon,x,y)|\leq C_j,\;\;\;\;\;\;0\leq x,y\leq 1,\epsilon\geq 0,j=0,1.
\end{eqnarray}
and $\psi^{(0,\nu)}(y,s)=\frac{\partial ^{\nu}\psi(y,s)}{\partial s^{\nu}}\leq\frac{ b^{\nu}}{E(\nu,\nu)},0\leq b\leq1,\;\;\;\;\;\nu=0,1,2,...$
where $E(n,k)\leq {{2n-1}\choose {n-1}}$ is the number of integer solutions of the equation
\begin{eqnarray}
&&\sum_{j=1}^{k}r_j\cdot s_j=n,0\leq r_j,s_j\leq n, j=1,2,..,k.\label{inteq0000}
\end{eqnarray}
Assuming that the unperturbed equation for $\epsilon=0$ has a solution in $C([0,1])$
then the perturbation series
%
%
\begin{eqnarray}
&&\phi_{(0,0)}(x)=\phi(0,x)=\omega \int_0^1\Gamma_0(x,y)\psi(y,\phi_{(0,0)}(y))dy\label{ppsiper0q}\\
&&\phi_{(0,1)}(x)=\omega\int_0^1\Gamma_1(x,y)\psi(y,\phi_{(0,0)}(y))dy\nonumber\\
&+&\omega\int_0^1\Gamma_0(x,y)\psi^{(0,1)}(y,\phi_{(0,0)}(y))\phi_{(0,1)}(y)dy\;\;\;\;\;\;\;\;\label{ppsiper1q}\\
&&\cdot\cdot\cdot\cdot\\
&&\phi_{(0,\nu)}(x)=\omega\int_0^1\Gamma_0(x,y)\psi^{(0,\nu)}(y,\phi_{(0,0)}(y))P_{\nu+1}(y)dy\nonumber\\
&+&\omega\int_0^1\Gamma_1(x,y)\psi^{(0,\nu-1)}(y,\phi_{(0,0)}(y))P_{\nu}(y)dy\;\;\;\;\;\;\;\;\;\;\;\label{ppsiper0nq}\\
&&P_{\nu}(y)=\nonumber\\
\sum_{{\sum_{j=1}^{\nu-1}r_j\cdot s_j=\nu-1}\atop{r_j,s_j=0,j=1,...,\nu-1}}^{\nu-1}&&{{\nu-1}\choose {r_1,s_1;\cdot\cdot\cdot;r_{\nu-1},s_{\nu-1}}}\prod_{j=1}^{\nu-1}\phi_{(0,s_j)}^{r_j}(y),\;\;\;\nu=2,3,...\;\;\;\;\;\label{auxppserq}\\
&&{m\choose {r_1,s_1;\cdot\cdot;r_m,s_m}}=\frac{m!}{\prod_{j=1}^m (s_j!)^{r_j}}\\
&&\psi^{(k,\nu)}(s,\phi(\epsilon,y))=\frac{\partial ^k\partial ^{\nu}\psi(s,\phi(\epsilon,y))}{\partial s^k\partial \epsilon^{\nu}},\;\;\;\;\;k,\nu=0,1,2,...\;\;\;\;\;\;\;\;\label{part2222}
\end{eqnarray}
%
%
converges absolutely  in      $||\cdot||_0=\sup_{x\in[0,1]}|\cdot(x)|$ norm and $||\cdot||_2$ norm to a solution
\begin{eqnarray}
&&\phi(\epsilon,x)=\sum_{j=0}^{\infty}\phi_{(0,j)}(x)\frac{\epsilon^j}{j!}\label{serper}
\end{eqnarray}
which is continuous with respect to $(\omega,\epsilon,x)$ for suitable ranges.}\label{pertheo2}\\\\}
{\bf Proof:}  As a result of the $L^2$ integrability we have that we have converging in $L^2$ sequences of kernels to $\Gamma_0,\Gamma_1$ respectively such that  $\Gamma_{0,n}(x,y)\leq C_n'\rightarrow C_{0,}\leq\infty$ and $\Gamma_{1,n}(x,y)\leq C_n'\rightarrow C_{1}\leq\infty$ and we have that $||\Gamma_{0,n}||_1\leq C_n, ||\Gamma_{1,n}||_1\leq C_n'$, which implies that inductively we can show that $||\phi_n(0,\cdot)||_0\leq\left[ \frac{\omega (C_0}{1-\omega C_0}+ \frac{\omega C_1}{1-\omega C_0}\right](n)! D^{n+1},n=1,2,...$, so for $\omega \leq g_{\pm}(C_0,C_1,D)$ we have the bound $||\phi_n(0,\cdot)||_0\leq n!D^n, n=1,2,...$,
where if $\Delta(C_0,C_1,D)=[2DC_0+C_0+C_1]^2-4(DC_1+2C_0^2)<0$ we have for any $\omega $ the bound and otherwise for $\Delta(C_0,C_1,D)>0$ we have for
\begin{eqnarray}
&&\omega<g_{\pm}(C_0,C_1,D)=\frac{2DC_0+C_0+C_1\pm\sqrt{\Delta(C_0,C_1,D)}}{2(DC_1+2C_0^2)}>0
\end{eqnarray}
The induction is based on the inequality
\begin{eqnarray}
&&|P_{\nu}(y)|\leq E(\nu-1,\nu-1)(\nu-1)!
\end{eqnarray}
and
\begin{eqnarray}
&&||\phi_{0,\nu}(\cdot)||_0\leq |\omega|\int_0^1|\Gamma_0(x,y)||\psi^{(0,\nu)}(y,\phi_{0,0}(y))||P_{\nu+1}(y)|dy\nonumber\\
&+&|\omega|\int_0^1|\Gamma_1(x,y)||\psi^{(0,\nu-1)}(y,\phi_{0,0)}(y))||P_{\nu}(y)|dy\\
&\leq&|\omega|C_0||\phi_{\nu}(0,\cdot)||_0\frac{ b^{\nu}}{E(\nu,\nu)}+|\omega|C_0|\frac{ b^{\nu}}{E(\nu,\nu)}E(\nu-1,\nu-1)\nu !D^{\nu}\nonumber\\
&+&|\omega|C_1||_0\frac{ b^{\nu}}{E(\nu,\nu)}(\nu-1)!D^{\nu-1}.
\end{eqnarray}
and  which leads to the convergent series
\begin{eqnarray}
&&||\phi_{(n,\epsilon)}(x)||_0=\sum_{j=0}^{\infty}||\phi_{n,0,j)}(x)||_0\frac{\epsilon^j}{j!}\\
&&\leq \sum_{j=0}^{\infty}k j!D^j\frac{\epsilon^j}{j!}||\phi_{0,n}||_0,n=0,1,2,...\label{serper0n2a}
\end{eqnarray}
converging in $||\cdot||_0$ for any value of $\omega\leq g_{\pm}(C_0,C_1,D),1/D>\epsilon\geq 0$.\\\\
So we have that the limit solution $||\phi(\epsilon,x)||_0<\infty$ is bounded  everywhere for $(\omega,x,\epsilon)$ in the allowed ranges.\\
As the convergence is uniform we have that $\phi(\epsilon,x)$ is continuous on the specified ranges.\\\\
In the same manner we have as a result that $\phi(\epsilon,x)$ belongs to $L_2[0,1]$.\\\\
{\zcorr{Let $f$, $\Gamma_0$ be continuous and $\Gamma_1$ to be $L^2$ integrable. Then $\Gamma (\epsilon,x,y)$ defines a continuous , bounded solution for $\epsilon<\epsilon_0>0$.}\label{corpert0}\\\\}
{\bf Proof:} If $\Gamma_0$ is continuous the solution $\phi_{(0)}(x)$ as in $L^2[0,1]$ it implies that $\phi_{(0)}(x)=f(x)+\omega\int_0^1\Gamma_0(x,y)\phi_{(0)}(y)dy$ as a sum of a continuous function and an integral with a continuous kernel is continuous. Applying theorem (\ref{pertheo2}) we have the sought implication.\\\\
{\zcorr{Let $f$ continuous and  $\Gamma_1$ to be $L^2$ integrable. Then $\Gamma (\epsilon,x,y)=\epsilon I+\epsilon(\Gamma_1-I)$ defines a continuous , bounded solution for $\epsilon<\epsilon_0>0$.}\label{corpert1}\\\\}
{\bf Proof:} Immediate from the above corollary (\ref{corpert0}).
Next we consider the behavior of the solution of the perturbed integral equation when the kernel has bounded uniformly derivative with respect to $x$. Then the perturbation series has a bounded derivative and uniformly convergent solution.\\\\
{\theo{Let's consider the perturbed integral equation
\begin{eqnarray}
&&  \phi(\epsilon,x)-\omega\int_0^1\Gamma (\epsilon,x,y)\phi(\epsilon,y)dy=f(x),\;\;\;\;x\in[0,1],\epsilon\geq 0\label{eqinte1nld}
                                      \end{eqnarray}
                                      where $\Gamma (\epsilon,x,y)=\Gamma_0(x,y)+\epsilon \Gamma_1(x,y)$ is continuous and $L^2$ integrable respectively functions of their variables and $\int_0^1\int_0^1\Gamma (\epsilon,x,y)dx dy<\infty$.
                                       Let $$\sup_{0\leq x\leq 1}\frac{\partial \Gamma_j(x,y)}{\partial x}\leq C_1>0,\;\;j=0,1.$$
                                       Suppose that
\begin{eqnarray}
 && |\Gamma_0(\epsilon,x,y)|\leq C,\;\;\;\;\;\;0\leq x,y\leq 1,\epsilon\geq 0.
\end{eqnarray}
Assuming that the unperturbed equation for $\epsilon=0$ has a solution in $C^1([0,1])$ or $L^2([0,1])$ respectively
then the perturbation series
%
%
\begin{eqnarray}
&&\phi_{(0,0)}(x)=\omega \int_0^1\Gamma_0,x,y)\phi_{(0,0)}(y)dy+f(x)\label{per0dd}\\
&&\phi_{(0,1)}(x)=\omega \int_0^1\Gamma_1(x,y)\phi_{(0,0)}(y)dy+\omega \int_0^1\Gamma_0(x,y)\phi_{(0,1)}(y)dy\label{per1dd}\\
&&\cdot\cdot\cdot\nonumber\\
&&\phi_{(0,\nu)}(x)=\frac{\partial^{\nu}\phi(x)}{\partial x^{\nu}}|_{x=0}=\omega \int_0^1\Gamma_1((x,y)\phi_{(0,\nu-1)}(y)dy\nonumber\\
&&+\omega \int_0^1\Gamma (0,x,y)\phi_{(0,\nu)}(y)dy,\nu=1,2,....\label{perndd}
\end{eqnarray}
converges absolutely  in  $||\cdot||_0=\sup_{x\in[0,1]}|\cdot(x)|$ norm and $||\cdot||_2$ norm to a solution
\begin{eqnarray}
&&\phi(\epsilon,x)=\sum_{j=0}^{\infty}\phi_{(0,j)}(x)\frac{\epsilon^j}{j!}\label{serperdd}
\end{eqnarray}
which is continuous with respect to $(\omega,\epsilon,x)$ and has a convergent series of derivatives.}\label{pertheo1dd}\\\\}
{\bf Proof:} As $\Gamma_0(x,y)\leq C$ and $\Gamma_1(x,y)\leq C$ we have that $||\Gamma_0||_1\leq C, ||\Gamma_1||_1\leq C$, which implies
\begin{eqnarray}
||\phi_{(0,\nu)}||_0&=&||\phi_{(0,0)}'(x)||_0\leq |\omega| ||\Gamma_1'||_1\cdot ||\phi_{(0,\nu-1)}||_0\nonumber\\
&+&|\omega| ||\Gamma_0'||_1\cdot ||\phi_{(0,\nu)}||_0\label{ineqper01dd}
\end{eqnarray}
leading to
\begin{eqnarray}
&&||\phi_{(0,\nu)}'||_0\leq\frac{|\omega|||\Gamma_1'||_1}{1-|\omega|||\Gamma_0'||_1}||\phi_{(0,\nu-1)}'||_2\label{ineqper02dd}
\end{eqnarray}
which implies that for $$\frac{|\omega|||\Gamma_1'||_1}{1-|\omega|||\Gamma_0'||_1}=\rho_0,$$
we have that
\begin{eqnarray}
||\phi'_{(\epsilon)}(x)||_0&=&\sum_{j=0}^{\infty}||\phi_{(0,j)}'(x)||_0\frac{\epsilon^j}{j!}\\
&&\leq \sum_{j=0}^{\infty}\rho_0^{j}\frac{\epsilon^j}{j!}||\phi_0'||_0\label{serper02dd}
\end{eqnarray}
the solution converges in $||\cdot||_0$ for any value of $\omega,\rho_0,\epsilon\geq 0$.\\\\
Also we observe that as $||\Gamma_j'||_2<\infty$ for $j=0,1$,
\begin{eqnarray}
||\phi_{(0,\nu)}'||_2&=&||\phi'_{(0,0)}(x)||_2\leq |\omega| ||\Gamma_1'||_2\cdot ||\phi_{(0,\nu-1)}'||_2\nonumber\\
&+&|\omega| ||\Gamma_0'||_2\cdot ||\phi_{(0,\nu)}'||_2\label{ineqper1dd}
\end{eqnarray}
leading to
\begin{eqnarray}
&&||\phi_{(0,\nu)}'||_2\leq\frac{|\omega|||\Gamma_1'||_2}{1-|\omega|||\Gamma_0'||_2}||\phi_{(0,\nu-1)}'||_2\label{ineqper2dd}
\end{eqnarray}
which implies that for $$\frac{|\omega|||\Gamma_1'||_2}{1-|\omega|||\Gamma_0'||_2}=\rho,$$
we have that
\begin{eqnarray}
&&||\phi'_{(\epsilon,0)}(x)||_2\leq\sum_{j=0}^{\infty}||\phi_{(0,\nu)}'(x)||_2\left[\frac{\epsilon^j}{j!}\right]\nonumber\\
&&\leq \sum_{j=0}^{\infty}|\rho^{j}\left[\frac{\epsilon^j}{j!}\right]||\phi_0'||_2\leq||\phi_0'||_2 e^{|\rho\epsilon|} \label{serper2dd}
\end{eqnarray}
the derivative of the solution converges for any value of $\omega,\rho,\epsilon\geq 0$.\\\\
As a convergent power series everywhere it is continuous with respect to $(x,\epsilon)$ as one can see immediately from uniform convergence.
Also we have that is continuous with respect to $\omega$.\\\\\\
Next we continue with the following theorem.\\\\
{\theo{Let's consider the perturbed integral equation
\begin{eqnarray}
&&  \phi(\epsilon,x)-\omega\int_0^1\Gamma (\epsilon,x,y)\phi(\epsilon,y)dy=f(x),\;\;\;\;x\in[0,1],\epsilon\geq 0\label{eqinte1nld2}
                                      \end{eqnarray}
                                      where $\Gamma (\epsilon,x,y)=\Gamma_0(x,y)+\epsilon \Gamma_1(x,y)$ is  $L^2$ integrable respectively functions of their variables and $\int_0^1\int_0^1\Gamma (\epsilon,x,y)dx dy<\infty$.
                                       Suppose that
\begin{eqnarray}
 && |\Gamma_j(x,y)|\leq C,\;\;\;\;\;\;0\leq x,y\leq 1,\epsilon\geq 0, j=0,1.
\end{eqnarray}
Assuming that the unperturbed equation for $\epsilon=0$ has a solution in $C^1([0,1])$
then the perturbation series
\begin{eqnarray}
&&\phi_{(0,0)}(x)=\omega \int_0^1\Gamma_0,x,y)\phi_{(0,0)}(y)dy+f(x)\label{per0d2}\\
&&\phi_{(0,1)}(x)=\omega \int_0^1\Gamma_1(x,y)\phi_{(0,0)}(y)dy+\omega \int_0^1\Gamma_0(x,y)\phi_{(0,1)}(y)dy\label{per1d2}\\
&&\cdot\cdot\cdot\nonumber\\
&&\phi_{(0,\nu)}(x)=\frac{\partial^{\nu}\phi(x)}{\partial x^{\nu}}|_{x=0}=\omega \int_0^1\Gamma_1((x,y)\phi_{(0,\nu-1)}(y)dy\nonumber\\
&&+\omega \int_0^1\Gamma (0,x,y)\phi_{(0,\nu)}(y)dy,\nu=1,2,....\label{pernd2}
\end{eqnarray}
converges absolutely  in      $||\cdot(x)||_0=\sup_{x\in[0,1]}|\cdot(x)|$ norm and $||\cdot||_2$ norm to a solution
\begin{eqnarray}
&&\phi(\epsilon,x)=\sum_{j=0}^{\infty}\phi_{(0,j)}(x)\frac{\epsilon^j}{j!}\label{serperd2}
\end{eqnarray}
which is continuous with respect to $(\omega,\epsilon,x)$.}\label{pertheo2d2}\\\\}
{\bf Proof:}  As  we have that $||\Gamma_{0}||_1\leq C, ||\Gamma_{1}||_1\leq C$, which implies corresponding inequalities to (\ref{ineqper01dd},\ref{ineqper02dd}), and therefore there exists a solution sequence
\begin{eqnarray}
&&||\phi_{(n,\epsilon)}(x)||_0=\sum_{j=0}^{\infty}||\phi_{(n,0,j)}(x)||_0\frac{\epsilon^j}{j!}\\
&&\leq \sum_{j=0}^{\infty}\rho_{0,n}^{j}\frac{\epsilon^j}{j!}||\phi_{0,n}||_0\leq ||\phi_{0,n}||_0e^{|\rho_{0,n}\epsilon},n=0,1,2,...\label{serper0n2d2}
\end{eqnarray}
the solution converges in $||\cdot||_0$ for any value of $\omega,\rho_0,\epsilon\geq 0$.\\\\
As $\rho_{0,n}\rightarrow \rho_0<\infty$ we have that the limit solution $||\phi(\epsilon,x)||_0<\infty$ is bounded  everywhere for $(\omega,x,\epsilon)$.
As the convergence is uniform we have that $\phi(\epsilon,x)$ is $C^1[0,1]$ continuous.\\\\
In the same manner we have as a result that $\phi'_{(\epsilon)}(x)$ belongs to $L_2[0,1]$.\\\\
We give the respective theorem to (\ref{pertheo2}) for kernels with uniformly bounded derivative with respect to $x$.\\\\
{\theo{Let's consider the perturbed integral equation
\begin{eqnarray}
&&  \phi(\epsilon,x)-\omega\int_0^1\Gamma (\epsilon,x,y)\psi(y,\phi(\epsilon,y))dy=f(x),\;\;\;\;x\in[0,1],\epsilon\geq 0\label{eqinte1nold}
                                      \end{eqnarray}
                                      where $\Gamma (\epsilon,x,y)=\Gamma_0(x,y)+\epsilon \Gamma_1(x,y)$ is  $L^2$ integrable respectively function of its variables and $\int_0^1\int_0^1|\Gamma (\epsilon,x,y)|^2dx dy<\infty$.
                                      Let $$\sup_{0\leq x\leq 1}\frac{\partial \Gamma_j(x,y)}{\partial x}\leq C_j'>0,\;\;j=0,1.$$
                                       Suppose that
\begin{eqnarray}
 && |\Gamma_j(x,y)|\leq C_j,\;\;\;\;\;\;0\leq x,y\leq 1,\epsilon\geq 0,j=0,1.\\
\end{eqnarray}
and
and $\psi^{(0,\nu)}(y,s)=\frac{\partial ^{\nu}\psi(y,s)}{\partial s^{\nu}}\leq\frac{ b^{\nu}}{E(\nu,\nu)},0\leq b\leq1,\;\;\;\;\;\nu=0,1,2,...$
where $E(n,k)\leq {{2n-1}\choose {n-1}}$ is the number of integer solutions of the equation
\begin{eqnarray}
&&\sum_{j=1}^{k}r_j\cdot s_j=n,0\leq r_j,s_j\leq n, j=1,2,..,k.\label{inteq0000}
\end{eqnarray}
Assuming that the unperturbed equation for $\epsilon=0$ has a solution in $C([0,1])$
then the perturbation series
\begin{eqnarray}
&&\phi_{(0,0)}(x)=\omega \int_0^1\Gamma_0(x,y)\psi(y,\phi_{(0,0)}(y))dy\label{ppsiper0qd}\\
&&\phi_{(0,1)}(x)=\omega \int_0^1\Gamma_1(x,y)\psi(y,\phi_{(0,0)}(y))dy\nonumber\\
&+&\omega\int_0^1\Gamma_0(x,y)\psi^{(0,1)}(y,\phi_{(0,0)}(y))\phi_{(0,1)}(y)dy\;\;\;\;\;\;\label{ppsiper1qd}\\
&&.....\nonumber\\
&&\phi_{(0,\nu)}(x)\nonumber\\
&&=\omega\int_0^1\Gamma_0(x,y)\psi^{(0,\nu)}(y,\phi_{(0,0)}(y))P_{\nu+1}(y)dy\nonumber\\
&+&\omega\int_0^1\Gamma_1(x,y)\psi^{(0,\nu-1)}(y,\phi_{(0,0)}(y))P_{\nu}(y)dy\;\;\;\;\;\;\label{ppsiper0nqda}\\
P_{\nu}(y)=\nonumber\\
\sum_{{\sum_{j=1}^{\nu-1}r_j\cdot s_j=\nu-1}\atop{r_j,s_j=0,j=1,...,\nu-1}}^{\nu-1}&&{{\nu-1}\choose {r_1,s_1;\cdot\cdot\cdot;r_{\nu-1},s_{\nu-1}}}\prod_{j=1}^{\nu-1}\phi_{(0,s_j)}(x)^{r_j}(y),\nu=2,3,...\;\;\label{auxppserqd}\\
&&{m\choose {r_1,s_1;\cdot\cdot;r_m,s_m}}=\frac{m!}{\prod_{j=1}^m (s_j!)^{r_j}}\\
&&\psi^{(k,\nu)}(s,\phi(\epsilon,y))=\frac{\partial ^k\partial ^{\nu}\psi(s,\phi(\epsilon,y))}{\partial s^k\partial \epsilon^{\nu}}k,\nu=0,1,2,...\;\;\;\;\;\;\;\;\label{part2222d}
\end{eqnarray}
%
%
converges absolutely  in      $||\cdot||_0=\sup_{x\in[0,1]}|\cdot(x)|$ norm and $||\cdot||_2$ norm to a solution
\begin{eqnarray}
&&\phi(\epsilon,x)=\sum_{j=0}^{\infty}\phi_{(0,j)}(x)\frac{\epsilon^j}{j!}\label{serperda0}
\end{eqnarray}
which is continuous with respect to $(\omega,\epsilon,x)$ for suitable ranges.}\label{pertheo2d}\\\\}
{\bf Proof:}  As $\Gamma_0(x,y)\leq C$ and $\Gamma_1(x,y)\leq C$ we have that $||\Gamma_0||_1\leq C, ||\Gamma_1||_1\leq C_1$, which implies that inductively we can show as in theorem (\ref{pertheo2}) that $||\phi_{(0,n)}(\cdot)||_0\leq\left[ \frac{\omega C_0}{1-\omega C_0}+ \frac{\omega C_1}{1-\omega C_0}\right](n)! D^{n},n=1,2,...$, so for $\omega \leq g_{\pm}(C_0,C_1,D)$ we have the bound $||\phi_{(0,n)}(\cdot)||_0\leq n!D^{n+1}, n=1,2,...$,
where if $\Delta(C_0,C_1,D)=[2DC_0+C_0+C_1]^2-4(DC_1+2C_0^2)<0$ we have for any $\omega $ the bound and otherwise for $\Delta(C_0,C_1,D)>0$ we have for
\begin{eqnarray}
&&\omega<g_{\pm}(C_0,C_1,D)=\frac{2DC_0+C_0+C_1\pm\sqrt{\Delta(C_0,C_1,D)}}{2(DC_1+2C_0^2)}>0
\end{eqnarray}
%
%
%
%
%
%
This leads to
\begin{eqnarray}
&&||\phi_{(\epsilon,n)}(x)||_0\leq\sum_{j=0}^{\infty}||\phi_{(0,j,n)}(x)||_0\frac{\epsilon^j}{j!}\\
&&\leq \sum_{j=0}^{\infty}k j!D^j\frac{\epsilon^j}{j!}||\phi_{0,n}||_0,n=0,1,2,...\label{serper0n2}
\end{eqnarray}
converging in $||\cdot||_0$ for any value of $\omega\leq g_{\pm}(C_0,C_1,D),1/D>\epsilon\geq 0$.\\\\
So we have that the limit solution $||\phi(\epsilon,x)||_0<\infty$ is bounded  everywhere for $(\omega,x,\epsilon)$ in the allowed ranges.\\
As the convergence is uniform we have that $\phi(\epsilon,x)$ is continuous on the specified ranges.\\\\
In the same manner we have as a result that $\phi(\epsilon,x)$ belongs to $L_2[0,1]$.\\
We also have in  the same manner that $\sup_{0\leq x\leq 1}|\frac{\partial \phi_{0,n}}{\partial x}|\leq n! D_1^n, n=2,3,...$
for $$\omega\leq g_{\pm}(C_0',C_1',D_1)$$ holds and  which leads to the uniformly convergent series
\begin{eqnarray}
&&||\phi_{(0,n)}'(x)||_0\leq\sum_{j=0}^{\infty}||\phi_{(0,j,n)}(x)||_0\frac{\epsilon^j}{j!}\nonumber\\
&&\leq \sum_{j=0}^{\infty}k_1 j!D_1^j\frac{\epsilon^j}{j!}||\phi'_{(0,n)}||_0,n=0,1,2,...\label{serper0n2d}
\end{eqnarray}
converging in $||\cdot||_0$ for any value of $\omega\leq g_{\pm}(C_0',C_1',D_1),1/D_1>\epsilon\geq 0$.\\\\
So we have that the limit solution $||\phi'(\epsilon,x)||_0<\infty$ is bounded  everywhere for $(\omega,x,\epsilon)$ in the allowed ranges.\\
As the convergence is uniform we have that $\phi'(\epsilon,x)$ is continuous on the specified ranges.\\\\
In the same manner we have as a result that $\phi'_{(\epsilon,0)}(x)$ belongs to $L_2[0,1]$.\\
\\\\
{\zcorr{Let $f$, $\Gamma_0$ be continuous and $\Gamma_1$ to be $L^2$ integrable.  Let $$\sup_{0\leq x\leq 1}\frac{\partial \Gamma_j(x,y)}{\partial x}\leq C_1>0,\;\;j=0,1.$$
Then $\Gamma (\epsilon,x,y)$ defines a continuous , bounded $C^1[0,1]$ solution for $\epsilon<\epsilon_0>0$.}\label{corpert0d}\\\\}
{\bf Proof:} If $\Gamma_0$ is continuous the solution $\phi(0,x)$ as in $L^2[0,1]$ it implies that $\phi(0,x)=f(x)+\omega\int_0^1\Gamma_0(x,y)\phi(0,y)dy$ as a sum of a continuous function and an integral with a continuous kernel is continuous. Applying theorem (\ref{pertheo2d}) we have the sought implication.\\\\
{\zcorr{Let $f$ continuous and  $\Gamma_1$ to be $L^2$ integrable.  Let $$\sup_{0\leq x\leq 1}\frac{\partial \Gamma_j(x,y)}{\partial x}\leq C_1>0,\;\;j=0,1.$$
Then $\Gamma (\epsilon,x,y)=\epsilon I+\epsilon(\Gamma_1-I)$ defines a continuous , $C^1[0,1]$ bounded solution for $\epsilon<\epsilon_0>0$.}\label{corpert1d}\\\\}
{\bf Proof:} Immediate from the above corollary (\ref{corpert0d}).
We have the respective theorem as well.\\\\
{\theo{Let's consider the perturbed integral equation
\begin{eqnarray}
&&  \phi(\epsilon,x)-\omega\int_0^1\Gamma (\epsilon,x,y)\psi(y,\phi(\epsilon,y))dy=f(x),\;\;\;\;x\in[0,1],\epsilon\geq 0\label{eqinte1nola}
                                      \end{eqnarray}
                                      where $\Gamma (\epsilon,x,y)=\Gamma_0(x,y)+\epsilon \Gamma_1(x,y)$ is  $L^1$ integrable respectively functions of their variables and $\int_0^1\int_0^1\Gamma (\epsilon,x,y)dx dy<\infty$.
                                       Let $$\sup_{0\leq x\leq 1}\frac{\partial \Gamma_j(x,y)}{\partial x}\leq C_2>0,\;\;j=0,1.$$
                                       Suppose that
\begin{eqnarray}
 && |\Gamma_j(\epsilon,x,y)|\leq C_j,\;\;\;\;\;\;0\leq x,y\leq 1,\epsilon\geq 0,j=0,1.
\end{eqnarray}
Assuming that the unperturbed equation for $\epsilon=0$ has a solution in $C([0,1])$
then the perturbation series
\begin{eqnarray}
&&\phi_{(0,0)}(x)=\phi(0,x)=\omega \int_0^1\Gamma_0(x,y)\psi(y,\phi_{(0,0)}(y))dy\label{ppsiper0qd2}\\
&&\phi_{(0,1)}(x)=\omega \int_0^1\Gamma_1(x,y)\psi(y,\phi_{(0,0)}(y))dy\nonumber\\
&+&\omega\int_0^1\Gamma_0(x,y)\psi^{(0,1)}(y,\phi_{(0,0)}(y))\phi_{(0,1)}(y)dy\;\;\;\;\;\;\;\;\label{ppsiper1qd2}\\
&&.....\nonumber\\
&&\phi_{(0,\nu)}(x)\nonumber\\
&&=\omega\int_0^1\Gamma_0(x,y)\psi^{(0,\nu)}(y,\phi_{(0,0)}(y))P_{\nu+1}(y)dy\nonumber\\
&+&\omega\int_0^1\Gamma_1(x,y)\psi^{(0,\nu-1)}(y,\phi_{(0,0)}(y))P_{\nu}(y)dy\;\;\;\;\;\;\;\;\;\label{ppsiper0nqd}\\
P_{\nu}(y)=\nonumber\\
\sum_{{\sum_{j=1}^{\nu-1}r_j\cdot s_j=\nu-1}\atop{r_j,s_j=0,j=1,...,\nu-1}}^{\nu-1}&&{{\nu-1}\choose {r_1,s_1;\cdot\cdot\cdot;r_{\nu-1},s_{\nu-1}}}\prod_{j=1}^{\nu-1}\phi_{(0,s_j)}(x)^{r_j}(y),\nu=2,3,...\;\;\;\label{auxppserqd2}\\
&&{m\choose {r_1,s_1;\cdot\cdot;r_m,s_m}}=\frac{m!}{\prod_{j=1}^m (s_j!)^{r_j}}\\
&&\psi^{(k,\nu)}(s,\phi(\epsilon,y))=\frac{\partial ^k\partial ^{\nu}\psi(s,\phi(\epsilon,y))}{\partial s^k\partial \epsilon^{\nu}},\;\;k,\nu=0,1,2,...\;\;\;\;\;\;\;\;\label{part2222d2}
\end{eqnarray}
%
%
converges absolutely and uniformly in      $||\cdot||_0=\sup_{x\in[0,1]}|\cdot(x)|$ norm and $||\cdot||_1$ norm to a  $C^1[0,1]$ solution
\begin{eqnarray}
&&\phi(\epsilon,x)=\sum_{j=0}^{\infty}\phi_{(0,j)}(x)\frac{\epsilon^j}{j!}\label{serperd}
\end{eqnarray}
which is continuous with respect to $(\omega,\epsilon,x)$ for suitable ranges.}\label{pertheo3d}\\\\}
{\bf Proof:} Along the same lines as in theorem (\ref{pertheo2d}).\\\\
{\bf Note:} We have corresponding results if   $$\sup_{0\leq x\leq 1}\frac{\partial^m \Gamma_j(x,y)}{\partial x^m}\leq C_m>0,\;\;j=0,1,m=1,2,...,p.$$
Then the respective nonlinear Fredholm equations have a $C^m[0,1]$ solution for $m=1,2,...,p.$
\section{Solution procedure for nonlinear integral equations}\label{sec3}
Here we investigate 'small' perturbations employed on the non linear term where the unknown function is involved rather than the kernel of the integral equation. Again we show in theorem (\ref{pertheo2ds2}) that the perturbative series converges to a continuous solution. Then we show in therem (\ref{pertheo2ds2m}) below next, the existence of a maximal perurbation range for the perurbation parameter using Zorn's lemma on a appropriate partially ordered set of partitions.\\\\
{\theo{Let's consider the nonlinear integral equation
\begin{eqnarray}
&&  \phi(\epsilon,x)-\omega\int_0^1\Gamma (x,y)\psi(\epsilon,y,\phi(\epsilon,y))dy=0,\;\;\;\;x\in[0,1],\epsilon\geq 0\label{eqinte1nolds}
                                      \end{eqnarray}
                                      where $\Gamma (x,y)$ is  $L^2$ integrable respectively function of its variables and\\
                                       $\int_0^1\int_0^1\Gamma (\epsilon,x,y)dx dy<\infty$.
                                      We have $\psi(\epsilon,y,z)=z+\epsilon\Psi(y,z)$.
                                       Suppose that
\begin{eqnarray}
 && |\Gamma (x,y)|\leq C,\;\;\;\;\;\;0\leq x,y\leq 1,\epsilon\geq 0.
\end{eqnarray}
and $\Psi^{(0,\nu)}(y,s)=\frac{\partial ^{\nu}\Psi(y,s)}{\partial s^{\nu}}\leq\frac{ b^{\nu}}{E(\nu,\nu)},0\leq b\leq1,\;\;\;\;\;\nu=0,1,2,...$
where $E(n,k)\leq {{2n-1}\choose {n-1}}$ is the number of integer solutions of the equation
\begin{eqnarray}
&&\sum_{j=1}^{k}r_j\cdot s_j=n,0\leq r_j,s_j\leq n, j=1,2,..,k.\label{inteq0000}
\end{eqnarray}
If the integral equation \begin{eqnarray}
&&\phi_0(0,x)=\phi(0,x)=\omega \int_0^1\Gamma (x,y)\phi_0(0,y)dy\label{ppsiper0000}\\
\end{eqnarray}
has continuous solutions in $C[0,1]$, then the integral equation (\ref{eqinte1nolds})  has a solution in $C[0,1]$.}\label{pertheo2ds2}\\\\}
{\bf Proof:} Differentiating eq.(\ref{eqinte1nolds}) with respect to $\epsilon$ at $0$, we get the series of integral equations
\begin{eqnarray}
&&\phi_{(0,0)}(x)=\phi(0,x)=\omega \int_0^1\Gamma (x,y)\phi_{(0,0)}(y)dy\label{ppsiper0}\\
&&\phi_{(0,1)}(x)=\omega \int_0^1\Gamma (x,y)\phi_{(0,1)}(y)dy\nonumber\\
&+&\omega\int_0^1\Gamma (x,y)\Psi(y,\phi_{(0,0)}(y))dy\;\;\;\;\;\;\label{ppsiper1}\\
&&.....\nonumber\\
&&\phi_{(0,\nu)}(x)=\omega\int_0^1\Gamma (x,y)\phi_{(0,\nu)}(y)dy\nonumber\\
&+&\omega\int_0^1\Gamma (x,y)\Psi^{(0,\nu-1)}(y,\phi_{(0,0)}(y))P_{\nu}(y)dy\;\;\;\;\;\;\label{ppsiper0n}\\
P_{\nu}(y)&=&\sum_{{\sum_{j=1}^{\nu-1}r_j\cdot s_j=\nu-1}\atop{r_j,s_j=0,j=1,...,\nu-1}}^{\nu-1}{{\nu-1}\choose {r_1,s_1;\cdot\cdot\cdot;r_{\nu-1},s_{\nu-1}}}\prod_{j=1}^{\nu-1}\phi_{(0,s_j)}^{r_j}(y),\nu=2,3,...\;\;\;\;\;\;\;\label{auxppser}\\
&&{m\choose {r_1,s_1;\cdot\cdot;r_m,s_m}}=\frac{m!}{\prod_{j=1}^m (s_j!)^{r_j}}
\end{eqnarray}
which are linear Fredholm integral equations of the second kind after the first.
The first has solutions the eigenfunctions $\phi_n,n=1,2,...$ of $\Gamma (x,y)$ with associated eigenvalues $\mu_n,n=1,2,...$ under suitable assumptions on $K$.
Equation (\ref{ppsiper0n}) has solution as a second kind Fredholm integral equation of the form
\begin{eqnarray}
&&\phi_n(0,x)=f_n(x)+\omega\sum_{j=1,\omega\mu_j\neq 1}\frac{\mu_j}{1-\omega\mu_j}(f_n,\phi_j)\phi_j(x)\label{seckfredsol}\\
&&f_n(x)=\omega\int_0^1\Gamma (x,y)\Psi^{(0,n-1)}(y,\phi_0(0,y))P_n(y)dy,n=1,2,...\label{auxfredsol}
\end{eqnarray}
Inductively we can show as in theorem (\ref{pertheo2}) that $||\phi_{(0,n)}(\cdot)||_0\leq \frac{\omega C}{1-\omega C}(n-1)! D^{n-1},n=1,2,...$,  so for $\omega \leq\frac{D}{C+CD}$ we have the bound $||\phi_{(0,n)}(\cdot)||_0\leq n!D^n, n=1,2,...$,
This implies that under the assumption for continuity of solution in eq.(\ref{ppsiper0000}) then the perturbation series solution converges absolutely and uniformly in      $||\cdot||_0=\sup_{x\in[0,1]}|\cdot(x)|$ norm and $||\cdot||_1$ norm to a  $C^1[0,1]$ solution
\begin{eqnarray}
&&\phi(\epsilon,x)=\sum_{j=0}^{\infty}\phi_{(0,j)}(x)\frac{\epsilon^j}{j!}\label{serperds2}
\end{eqnarray}
which is continuous with respect to $(\omega,\epsilon,x)$ for suitable ranges. More precisely it holds for $\epsilon<\frac{1}{D}$ and then the solution is absolutely bounded by $\frac{1}{1-D\epsilon}$ for all $x\in [0,1]$.\\\\
The question arises that as given a perturbative series solution to a non linear integral equation if the region of perturbation can be expanded further so as to extend the region of convergence of the series.\\
The following theorem establishes the existence of a maximal perturbation along a univariate family of perturbation functions.
{\theo{Let's consider the nonlinear integral equation
\begin{eqnarray}
&&  \phi(\epsilon,x)-\omega\int_0^1\Gamma (x,y)\psi(\epsilon,y,\phi(\epsilon,y))dy=0,\;\;\;\;x\in[0,1],\epsilon\geq 0\label{eqinte1noldsm}
                                      \end{eqnarray}
                                      where $\Gamma (x,y)$ is  $L^2$ integrable respectively functions of their variables and\\ $\int_0^1\int_0^1\Gamma (\epsilon,x,y)dx dy<\infty$.
                                      We have $\psi(\epsilon,y,z)=z+\epsilon\Psi(y,z)$.
                                       Suppose that
\begin{eqnarray}
 && |\Gamma_0(\epsilon,x,y)|\leq C,\;\;\;\;\;\;0\leq x,y\leq 1,\epsilon\geq 0.
\end{eqnarray}
and $\Psi^{(0,\nu)}(y,s)=\frac{\partial ^{\nu}\Psi(y,s)}{\partial s^{\nu}}\leq\frac{ b^{\nu}}{E(\nu,\nu)},0\leq b\leq1,\;\;\;\;\;\nu=0,1,2,...$
where $E(n,k)\leq {{2n-1}\choose {n-1}}$ is the number of integer solutions of the equation
\begin{eqnarray}
&&\sum_{j=1}^{k}r_j\cdot s_j=n,0\leq r_j,s_j\leq n, j=1,2,..,k.\label{inteq0000}
\end{eqnarray}
If the integral equation \begin{eqnarray}
&&\phi_0(0,x)=\phi(0,x)=\omega \int_0^1\Gamma (x,y)\phi_0(0,y)dy\label{ppsiper0000m}\\
\end{eqnarray}
has continuous solutions in $C[0,1]$, then the integral equation (\ref{eqinte1noldsm})  has a maximal perturbation solution in $C[0,1]$.}\label{pertheo2ds2m}\\\\}
{\bf Proof:} Let $\frac{\partial^j\psi(\epsilon,y,\phi(\epsilon,y))}{\partial \epsilon^j}\leq B(j,\epsilon))$ where $B(j,\epsilon)\leq Q b(\epsilon)^j)/E(j,j),Q>0,j=1,2,...$.
Let a convex, respectively concave function $f_1,f_2:[a,b]\rightarrow \mathbf{R}$ on its domain of definition. Then we define the functional on its domain
\begin{eqnarray}
  V_1(f_1,P,a,b) &=& \sum_{j=0}^{n-1}\left|\frac{f_1(x_{j+1})-f_1(x_j)}{x_{j+1}-x_j}\right|\label{discvar0} \\
  V_2(f_2,P,a,b)&=&\sum_{j=1}^n\left|\frac{f_2(x_{j+1})-f_2(x_j)}{x_{j+1}-x_j}\right|\label{discvar1} \\
  P &=& \{a=x_0,x_1,...,x_k,..x_n=b\}\label{part}
\end{eqnarray}
We observe that as $||P||=\max_{0\leq j\leq n-1}|x_{j+1}-x_j|\rightarrow 0$ as $n\rightarrow \infty$ we have that $\lim_{n\rightarrow \infty}V_j(f_j,P,a,b)=\int_a^b|f_j'(x)|dx,\;\;j=1,2$.
So constructing the functional $V(f,a,b)$  for a function $f$ with convex or concave on successive subintervals on its domain of definition $[a,b]$ as
\begin{eqnarray}
&&V(f,P,a,b)=\sum_{j=0}^{k}V_{i_j}(f_{i_j},P_j,a_j,b_j)\label{partsumconv}\\
&&\bigcup_{j=0}^k[a_j,b_j)=[a,b)\\
&&f_{i_j}\;\;\;\mbox{convex on} [a_j,b_j), \mod(j,2)=0\\
&&f_{i_j}\;\;\;\mbox{concave on} [a_j,b_j), \mod(j,2)=1\\
&&\mbox{ or reversely},\;\;\;j=0,1,...,k.\\
&&P=\bigcup_{j=0}^kP_j,\mbox{span}(P_j)=b_j-a_j,j=0,...,k.
\end{eqnarray}
We have that $\lim_{n\rightarrow \infty}V(f,P,a,b)=\int_a^b|f'(x)|dx,\;\;$.
Moreover from the construction of the discrete functionals the convergence is monotone decreasing. Additionally this allows monotonicity and convergence with respect to refinement in partitions of the interval $[a,b)$ and shall be used in the following construction.\\\\
Given a function $f:[0,\infty)\rightarrow \mathbf{R}$ we impose an order on the partitions of $\mathbf{R}^+$ considering $P\preceq Q$ if $\mbox{span}(P)< \mbox{span}(Q)$ or if $\mbox{span}(P)= \mbox{span}(Q)$ and $V(f,P)>V(f,Q)$.\\\\
Thus as each chain of ascending partitions $P_0\preceq P_1\preceq\cdot\cdot\cdot\preceq P_m$ has an upper bound with span $\sup_{j=0}^{\infty}\mbox{span}(P_j)$ or if spans are finally equal, with value of the functional $\liminf_{j\rightarrow\infty}V(f,P_j)$, by Zorn's lemma the partial ordered set of chains has a maximal element $P^*$. So we have a maximal perturbation range $\mbox{span}(P^*)$.
Suppose we are achieving a perturbation of magnitude $\epsilon_0>0$. And from this value we achieve an additional perturbation of magnitude $\epsilon_1-\epsilon_0$ with nonlinear term $\psi(\epsilon_1,y,\phi(\epsilon_1,y))$. Continuing in this manner we achieve  perturbations with magnitudes $\epsilon_0<\epsilon_1<\cdot\cdot\cdot<\epsilon_m,m\in\mathbf{N}$. This defines a partition $P=\{\epsilon_0,\epsilon_1,\cdot\cdot\cdot,\epsilon_m\},m\in\mathbf{N}$of $[0,\infty)$. We order in the same way as above with $\mbox{span}(P_m)$ and $V(\psi(\epsilon,\cot,\cdot),P)$ these partitions defined by successive perturbations of the integral equation (\ref{eqinte1noldsm}) and we have as a result that there is a maximal element by Zorn's lemma on this set of partitions $P_{\psi}$ with maximal span and therefore range.
\section{Conclusions}
We investigate perturbative solutions for kernel perturbed integral equations and prove the convergence in an appropriate ranges of the perturbation series. Next we investigate perturbation series solutions for nonlinear perturbations of integral equations of Hammerstein type and formulate conditions for their convergence.
Then applying Zorn's lemma to a construction for the succesive 'small' perturbations to the perturbing function $\psi(\epsilon,y,z)$ with respect to the variable $\epsilon$ we get that  there is a maximal span where $\epsilon$ can vary and give finite solutions to the integral equation (\ref{eqinte1noldsm}) which are useful for the practitioner. Furthermore we note that these methods can be extended to more general equations arising in mathematical physics and mechanics.\\\\

\end{document}